\documentclass[reqno,11pt]{amsart}

\usepackage{amssymb,dsfont,enumerate}
\usepackage[pagebackref]{hyperref}
\usepackage[portrait,a4paper,margin=3cm]{geometry}

\usepackage{graphicx,amssymb,amsmath}
\usepackage[all]{xy}

\newtheorem{theorem}{Theorem}%[section]

\newtheorem{lemma}[theorem]{Lemma}
\newtheorem{proposition}[theorem]{Proposition}

\newcommand{\dE}{\mathbb{E}}

\newcommand{\dP}{\mathbb{P}}
\newcommand{\dR}{\mathbb{R}}

\newcommand{\dZ}{\mathbb{Z}}
\newcommand{\cE}{\mathcal{E}}

\newcommand{\cM}{\mathcal{M}}

\newcommand{\ABS}[1]{{{\left| #1 \right|}}} % |1|
\newcommand{\BRA}[1]{{{\left\{#1\right\}}}} % {1}
\newcommand{\PAR}[1]{{{\left(#1\right)}}} % (1)
\newcommand{\veps}{\varepsilon}

\renewcommand{\leq}{\leqslant}             % redef. of < or =
\renewcommand{\geq}{\geqslant}             % redef. of > or =
\newcommand{\IND}{\mathbf{1}}

\newcommand{\TR}{\mathrm{tr}}

\newcommand{\rank}{\mathrm{rank}}

\begin{document}

\title[Euclidean random matrices]{On Euclidean random matrices in high dimension}
\author{Charles Bordenave}
\email{charles.bordenave(at)math.univ-toulouse.fr}
%\urladdr{http://www.math.univ-toulouse.fr/~bordenave/}
\address{CNRS \& Universit\'e de Toulouse, Institut de Math\'ematiques de Toulouse, 118 route de Narbonne, 31062 Toulouse, France}

\keywords{Euclidean random matrices, Marcenko-Pastur distribution, Log-concave distribution.}

\subjclass[2000]{60B20 ; 15A18.}

\maketitle

\begin{abstract}
In this note, we study the $n \times n$ random Euclidean matrix whose entry $(i,j)$ is equal to $f ( \| X_i - X_j \|  )$ for some function $f$ and the $X_i$'s are i.i.d.~isotropic vectors in $\dR^p$. In the regime where $n$ and $p$ both grow to infinity and are proportional, we give some sufficient conditions for the empirical distribution of the eigenvalues to converge weakly. We illustrate our result on log-concave random vectors. 
\end{abstract}

\section{Introduction}

Let $Y$ be an \emph{isotropic} random vector in $\dR^p$,  i.e.  $\dE Y = 0$, $\dE [ Y Y^T ]  = I / p $, where $I$ is the identity matrix. Let $(X_1, \cdots, X_n)$  be independent copies of $Y$. We define the $n \times n$ matrix $A$ by, for all $1 \leq i , j \leq n$, 
$$
A_{i j } =f ( \| X_i - X_j \| ^2 ), 
$$
where $f : [ 0, \infty) \to \dR$ is a measurable function and $\| \cdot \|$ denotes the Euclidean norm. The matrix $A$ is a random Euclidean matrix. It has already attracted some attention see e.g.~M\'ezard, Parisi and Zhee \cite{MR1724455}, Vershik \cite{MR2086637} or Bordenave \cite{MR2462254} and references therein.

If $B$ is a  symmetric matrix of size $n$, then its eigenvalues, say $\lambda_1(B), \cdots, \lambda_n(B)$ are real. The empirical spectral distribution (ESD) of $B$ is classically defined as 
$$
\mu_B = \frac 1 n \sum_{i=1} ^ n \delta_{\lambda_i(B)},
$$
where $\delta_x$ is the Dirac delta function at $x$. In this note, we are interested in the asymptotic convergence of $\mu_A$ as $p$ and $n$ converge to $+\infty$. This regime has notably been previously considered in El Karoui \cite{MR2589315} and Do and Vu \cite{DoVu}. More precisely, we fix a sequence $p(n)$ such that 
\begin{eqnarray}\label{eq:povern}
\lim_{n \to \infty} \frac{ p(n) }{ n }  =  y \in (0, \infty).
\end{eqnarray}
Throughout this note, we consider, on a common probability space, an array of random variables $(X_k(n))_{ 1 \leq k \leq n}$ such that $(X_1(n), \cdots, X_n(n))$ are independent copies of $Y(n)$, an isotropic vector in $\dR^{p(n)}$. For each $n$, we define the Euclidean matrix $A(n)$ associated. For ease of notation, we will often remove the explicit dependence in $n$: we write $p$, $Y$, $X_k$ or $A$ in place of $p(n)$, $Y(n)$, $X_k (n)$ or $A(n)$.

The Marcenko-Pastur probability distribution with parameter $1/y$ is given by 
$$
\nu_{MP}  (dx) =  (1 - y )^+ \delta_0 (dx) +  \frac{ y } { 2 \pi x } \sqrt { (y_+ - x) ( x - y_-) } \IND_{ [ y_- , y_+] } (x) dx, 
$$
where $x^+ = ( x \vee 0)$, $y_ {\pm} = ( 1  \pm \frac 1  {\sqrt y }) ^2$ and $dx$ denotes the Lebesgue measure. Since the celebrated paper  of Marcenko and Pastur \cite{MR0208649}, this distribution is known to be closely related to empirical covariance matrices in high-dimension.

We say that $Y$ has a \emph{log-concave distribution}, if $Y$ has a density on $\dR^p$ which is log-concave. Log-concave random vectors have an increasing importance in convex geometry, probability and statistics (see e.g. Barthe \cite{MR2648682}). We will prove the following result. 

\begin{theorem}\label{th:main}
If $Y$ has a log-concave distribution and $f$ is three times differentiable at $2$, then, almost surely, as $n \to \infty$, 
$
\mu_{A}$ converges weakly to $\mu$, the law of $f(0 ) - f(2) +2 f'(2) - 2 f'(2) S$, where $S$ has distribution $\nu_{MP}$. 
\end{theorem}

With the weaker assumption that $f$ is differentiable at $2$, Theorem \ref{th:main} is conjectured in Do and Vu \cite{DoVu}. Their conjecture has motivated this note. It would follow from the thin-shell hypothesis which asserts that there exists $c > 0$, such that for any isotropic log-concave vector $Y$ in $\dR^p$, $\dE (\| Y\| -1 )^2 \leq c/p$ (see Anttila, Ball and Perissinaki \cite{MR1997580} and Bobkov and Koldobsky \cite{MR2083387}).  Klartag \cite{MR2520120} has proved the thin-shell hypothesis for isotropic unconditional log-concave vectors.

The proof of Theorem \ref{th:main} will rely on two recent results on log-concave vectors. Let $X = X(n) $ be the $n \times n$ matrix with columns given by $(X_1(n), \cdots, X_n (n))$. Pajor and Pastur have proved the following : 
\begin{theorem}[\cite{MR2539559}] \label{th:PP} If $Y$ has a log-concave distribution, then, in probability, as $n \to \infty$,
$ \mu_{ X  ^T X }$ converges weakly to  $ \nu_{MP}$. 
\end{theorem}

We will also rely on a theorem due to Gu\'edon and Millman.

\begin{theorem}[\cite{MR2846382}]\label{th:GM}
There exist positive constants $c_0, c_1$ such that if $Y$ is an isotropic log-concave vector in $\dR^p$, for any $t \geq 0$,
$$
\dP \PAR{ \ABS{  \| Y \| - 1 } \geq t  } \leq c_1 \exp \PAR{  - c_0 \sqrt p \PAR{ t \wedge t^3 } }. 
$$
\end{theorem}

With Theorems \ref{th:PP} and \ref{th:GM} in hand, the heuristic behind Theorem \ref{th:main} is simple. Theorem \ref{th:GM} implies that $\|X_i \|^2 \simeq 1$ with high probability. Hence, since $\| X_i  - X_j \|^2  = \|X_i \|^2 + \|X_j \|^2 - 2 X_i ^T X_j $, a Taylor expansion of $f$ around $2$ gives 
$$
A_{i j } \simeq \left\{ \begin{array}{ll}
f ( 2 ) - 2 f' (  2) X_i ^T X_j  & \hbox{ if }  i \ne j \\
f(0) & \hbox{ if }  i = j.
\end{array} \right.
$$
In other words, the matrix $A$ is close to the matrix
\begin{equation}\label{eq:defM}
M  = ( f(0) - f(2) + 2 f'(2)  ) I + f ( 2 ) J  - 2 f'(2) X^T X, 
\end{equation}
where $I$ is the identity matrix and $J$ is the matrix with all entries equal to $1$. From Theorem \ref{th:PP}, $\mu_{X^TX}$ converges weakly to $\nu_{MP}$. Moreover, since $J $ has rank one, it is negligible for the weak convergence of ESD. It follows that $\mu_{M}$ is close to $\mu$. The actual proof of Theorem \ref{th:main} will be elementary and it will follow this heuristic. We shall use some standard perturbation inequalities for the eigenvalues. The idea to perform a Taylor expansion was already central in \cite{MR2589315,DoVu}.

Beyond Theorems \ref{th:PP}-\ref{th:GM}, the proof of Theorem \ref{th:main} is not related to log-concave vectors. In fact, it is nearly always  possible to linearize $f$ as soon as the norms of the vectors concentrate around their mean. More precisely, let us say that two sequences of probability measures $(\mu_n)$, $(\nu_n)$, are asymptotically weakly equal, if for any bounded continuous function $f$, $\int f d \mu_n - \int f d \nu_n$ converges to $0$. 

\begin{theorem}\label{th:main2}
Assume that there exists an integer $\ell \geq 1$ such that $\dE \ABS{ \| Y \| -1 }^{2 \ell} = O ( p^{-1} )$, and that for any $\veps > 0$, 
\begin{equation}\label{eq:boundEE}
\lim_{n \to \infty} \dP \PAR{ \max_{1 \leq i , j \leq n  } \BRA{ \ABS{ \| X_i - X_j \|^2 -  2 } \vee \ABS{ \| X_i\| ^2 - 1} } \leq \veps } = 1.
\end{equation}
Then, if $f$ is $\ell$ times differentiable at $2$,  almost surely,
$
\mu_{A}$ is asymptotically weakly equal to the law of $f(0 ) - f(2) +2 f'(2) - 2 f'(2) S$, where $S$ has distribution $\dE \mu_{X^T X}$. 
\end{theorem}
The case $\ell =1$ of the above statement is contained in Do and Vu \cite[Theorem 5]{DoVu}. Besides Theorem \ref{th:PP}, some general conditions on the matrix $X$ guarantee the convergence of $\mu_{X^T X}$, see Yin and Krishnaiah \cite{MR816299}, G{\"o}tze and Tikhomirov \cite{MR2092202} or Adamczak \cite{MR2820070}.

\section{Proofs}

\subsection{Perturbation inequalities}

We first recall some basic perturbation inequalities of eigenvalues and introduce a good notion of distances for ESD.  For $\mu$, $\nu$ two real probability measures, the \emph{Kolmogorov-Smirnov distance} can be defined as
\begin{equation*}\label{eq:KSdual}
d_{KS} ( \mu , \nu)    = \sup \left\{ \int f d \mu   - \int f d \nu  : \| f \|_{BV} \leq 1 \right\}, 
\end{equation*}
where, for $f : \dR \to \dR$, the bounded variation norm is 
$
\| f \|_{BV}  = \sup \sum_{ k \in \dZ} | f ( x_{k+1}) - f (x_k) |,  
$
and the supremum is over all real increasing sequence $(x_k)_{ k \in \dZ}$. 
The following inequality is a classical consequence of the interlacing of eigenvalues (see e.g.  Bai and Silverstein \cite[Theorem A.43]{MR2567175}). 
\begin{lemma}[Rank inequality]\label{le:rank}
If $B$, $C$ are $n \times n$ Hermitian matrices, then, 
$$
d_{KS} (\mu_B, \mu_C) \leq \frac{ \rank (  B - C )  }n.
$$
\end{lemma}

For $p \geq 1$, let $\mu$, $\nu$ be two real probability measures such that $\int |x|^p d \mu$ and  $\int |x|^p d \nu$  are finite. We define the \emph{$L^p$-Wasserstein distance} as
\begin{equation*}\label{eq:defWp}
W_p ( \mu, \nu)  = \PAR{ \inf_\pi    \int_{\dR \times \dR} | x - y| ^p  d \pi   } ^{\frac 1 p} 
\end{equation*}
where the infimum is over all coupling $\pi$ of $\mu$ and $\nu$ (i.e. $\pi$ is probability measure on $\dR \times \dR$ whose first marginal is equal to $\mu$ and second marginal is equal to $\nu$). H\"older inequality implies that for $1 \leq p \leq q$,
$W_p \leq W_{q}$. 
Moreover, the Kantorovich-Rubinstein duality gives a variational expression for $W_1$:
\begin{equation*}\label{eq:KRdual}
W_1 (\mu, \nu) =  \sup \left\{ \int f d \mu   - \int f d \nu  : \| f \|_{L} \leq 1 \right\},
\end{equation*}
where $\| f \|_L = \sup_{ x \ne y } |f(x) - f(y) | / |x- y| $ is the Lipschitz constant of $f$. The next classical inequality is particularly useful (see e.g. Anderson, Guionnet and Zeitouni \cite[Lemma 2.1.19]{AGZ}). 
\begin{lemma}[Hoeffman-Wielandt inequality]\label{le:HW}
If $B$, $C$ are $n \times n$ Hermitian matrices, then 
$$
W_2 ( \mu_B, \mu_C) \leq \sqrt{ \frac 1 n \TR ( B - C)^2}. 
$$
\end{lemma}

We finally introduce the distance
$$
d ( \mu, \nu) = \sup \BRA{ \int f d \mu - \int f d \nu :  \| f \|_L \leq 1  \hbox{ and } \| f \|_{BV} \leq 1}. 
$$
By Lemmas \ref{le:rank} and \ref{le:HW}, we obtain that for any $n \times n$ Hermitian matrices $B$, $C$,
\begin{equation}\label{eq:hehe}
d ( \mu_B , \mu _C) \leq \sqrt{ \frac 1 n \TR ( B - C)^2} \wedge \frac{ \rank ( B- C )  }n.
\end{equation}
Notice that $d( \mu_n , \mu) \to 0$ implies that $\mu_n$ converges weakly to $\mu$.

\subsection{Concentration inequality}

For $x  = ( x_1, \cdots, x_n  )  \in \cM_{p , n} (\dR)$, define $a(x)$ as the Euclidean matrix obtained from the columns of  $x$ : $a(x) _{ij} = f ( \| x_i - x_j \|^2)$. In particular, we have $A = a (X)$. Let $i \in \{1, \cdots, n\}$, $x' = (x'_1, \cdots, x'_n) \in \cM_{p , n} (\dR)$ and assume that $x'_j = x_j$ for all $j \ne i$. Then $a(x)$ and $a(x')$ have all entries equal but the entries on the $i$-th row or column. We get
$$
\rank ( a(x) - a(x') ) \leq 2. 
$$
It thus follows from Lemma \ref{le:rank} that for any function $f$ with $\| f \|_{BV} < \infty$, 
$$
\left| \int f d \mu_{a(x) }  - \int f d \mu_{a(x') } \right| \leq \frac { 2 \| f \|_{BV} } { n}. 
$$
Using Azuma-Hoeffding's inequality, it is then straightforward to check that for any $t \geq 0$, 
\begin{equation}\label{eq:concineq}
\dP \PAR{ \int f d \mu_{A }  - \dE  \int f d \mu_{A }  \geq t } \leq \exp\PAR{  - \frac {  n t ^2 }{ 8 \| f \|^2_{BV} } }. 
\end{equation}
(For a proof, see \cite[proof of Lemma C.2]{MR2837123} or Guntuboyina and Leeb \cite{MR2535081}). Using the Borel-Cantelli Lemma, this shows that for any such function $f$, a.s. 
$$
\int f d \mu_{A }  -    \int f d \dE \mu_{A } \to 0. 
$$

Now, recall that $M$ was defined by \eqref{eq:defM}. Since the matrix $J$ has rank one, from Theorem \ref{th:PP} and Lemma \ref{le:rank}, $\dE \mu_M$ converges weakly to $\mu$.  Hence our Theorem \ref{th:main} is a corollary of the following proposition. 

\begin{proposition}\label{prop:main}
Under the assumptions of Theorem \ref{th:main}, we have
$$
\lim_{n \to \infty} d \PAR{ \dE \mu_{A} , \dE \mu_{M} } = 0. 
$$
\end{proposition}

\subsection{Proof of Proposition \ref{prop:main}} The idea is to perform a multiple Taylor expansion which takes the best out of \eqref{eq:hehe}. 

 \subsubsection*{Step 1 : concentration of norms} 

By assumption, there exists an open interval $K  = (2 - \delta , 2 + \delta)$ such that $f$ is $C^1$ in $K$ and, for any $x \in K$, 
$$
f( x ) = f(2) + f'(2) (x-2) + \frac{f'' (2)}{2} ( x - 2)^2 + \frac{f''' (2)}{6} ( x - 2)^3 ( 1 + o ( 1) ).   
$$

For any $i \ne j$, $(X_i- X_j)/\sqrt 2$ is an isotropic log-concave vector. Define the sequence $\veps (n) = n ^ { - \kappa } \wedge ( \delta /2)$ with $0 < \kappa < 1/6$. It follows from Theorem \ref{th:GM} and the union bound that the event 
$$
\cE = \BRA{ \max_{i , j } \BRA{ \ABS{ \| X_i - X_j \|^2 - 2 } \vee \ABS{ \| X_i\|^2  - 1} } \leq \veps(n) }
$$
has probability tending to $1$ as $n$ goes to infinity.

\subsubsection*{Step 2 : Taylor expansion around $\| X_i \|^2 + \| X_j \|^2$} 

We consider the matrix 
$$
B_{i j } = \left\{ \begin{array}{ll}
f ( \| X_i \|^2 + \|X_j\|^2 ) - 2 f' (  \| X_i \|^2 + \|X_j\|^2 ) X_i ^T X_j  & \hbox{ if }  i \ne j \\
f(0) & \hbox{ if }  i = j.
\end{array} \right.
$$
On the event $\cE$, $\| X_i \|^2 + \|X_j\|^2 \in K$. Since $f$ is $C^1$ in $K$, we may perform a Taylor expansion of $f ( \| X_i - X_j \|^2)$ around $\|X_i \|^2 + \| X_j \|^2$. It follows that for $i \ne j$, 
$$
\ABS{ A_{ij} - B_{ij} } = o  \PAR{ \| X_i - X_j \|^2 - \|X_i \|^2 - \|X_j \|^2 }  \leq  \delta(n) \ABS{  X_i ^T X_j }, 
$$
where $\delta(n)$ is a sequence going to $0$.
From \eqref{eq:hehe} and Jensen's inequality, we get 
\begin{eqnarray*}
d ( \dE \mu_A , \dE \mu_B )  \leq  \dE d (  \mu_A ,  \mu_B )  & \leq & \dP ( \cE^c ) +\PAR{  \frac 1 n \sum_{ i \ne j } \dE |A_{ij} - B_{ij} |^2 \IND_\cE }^{1/2} \\
& \leq & \dP ( \cE^c ) +  \delta(n) \PAR{  n   \dE \ABS{  X_1 ^T X_2 }^2 }^{1/2}.
\end{eqnarray*}
Now, from the assumption that $X_1$ and $X_2$ are independent and isotropic, we find
\begin{eqnarray*}
\dE \ABS{  X_1 ^T X_2 }^2 & = & \dE \PAR{  \sum_{k=1} ^ p X_{k 1} X_ {k 2} } ^2 
 =     \sum_{k=1} ^ p \PAR { \dE X_{k 1}^2  } ^2 
 =   \frac 1 p.
\end{eqnarray*}
By assumption \eqref{eq:povern}, we deduce that 
$$
\lim_{ n \to \infty} d ( \dE \mu_A , \dE \mu_B )  = 0. 
$$
It thus remains to compare $\dE \mu_B$ and $\dE \mu_{M}$.

\subsubsection*{Step 3 : Taylor expansion around $2$}
We define the matrix 
$$
C_{i j } = \left\{ \begin{array}{ll}
f ( \| X_i \|^2 + \|X_j\|^2 ) - 2 f' (  2 ) X_i ^T X_j  & \hbox{ if }  i \ne j \\
f(0) & \hbox{ if }  i = j.
\end{array} \right.
$$
We now use the fact that $f'$ is locally Lipschitz at 2. It follows that if $\cE$ holds, for $i \ne j$, 
$$
\ABS{ B_{ij} - C_{ij} } = O  \PAR{  X_i^ T X_j (  \|X_i \|^2 + \|X_j \|^2 -2  ) }  \leq c\, \veps(n) \ABS{  X_i ^T X_j }. 
$$
The argument of step 2 implies that 
$$
\lim_{ n \to \infty} d ( \dE \mu_B , \dE \mu_C )  = 0. 
$$
It thus remains to compare $\dE \mu_C$ and $\dE \mu_{M}$. 

\subsubsection*{Step 4 : Taylor expansion around $2$ again}
We now consider the matrix
$$
D_{i j } = \left\{ \begin{array}{ll}
f ( 2 ) + f'(2) ( \| X_i \| ^2 + \| X_j \| ^2  -2 ) + \frac{ f''(2) } {2}  ( \| X_i \| ^2 + \| X_j \| ^2  -2 )^2 & \\ 
\quad \quad +\frac{ f'''(2) } {6}  ( \| X_i \| ^2 + \| X_j \| ^2  -2 )^3 -   2 f' (  2 ) X_i ^T X_j  & \hbox{ if }  i \ne j \\
f(0) & \hbox{ if }  i = j.
\end{array} \right.
$$
We are going to prove that 
\begin{equation}\label{eq:T3}
\lim_{ n \to \infty} d ( \dE \mu_C , \dE \mu_D )  = 0. 
\end{equation}

We perform a Taylor expansion of order $3$ of $f( \| X_i \| ^2+ \| X_j \|^2)$ around $2$. It follows that if $\cE$ holds, for $i \ne j$, 
$$
\ABS{ C_{ij} - D_{ij} } = o  \PAR{   \|X_i \|^2 + \|X_j \|^2 -2 }^3   \leq \delta (n) \ABS{   \|X_i \|^2 + \|X_j \|^2 -2 }^3,
$$
where $\delta(n)$ is a sequence going to $0$.  Using \eqref{eq:hehe} and arguing as in step 2, in order to prove \eqref{eq:T3}, it thus suffices to show that 
$$
\frac 1 n \sum_{ i \ne j } \dE | \|X_i \|^2 + \|X_j \|^2 -2  |^6 \IND_\cE= O (1). 
$$
Since, for $\ell \geq 1$, $| x + y | ^\ell \leq 2^{\ell-1} ( |x|^\ell + |y|^\ell)$, it is sufficient to show that 
$$
n \dE\PAR{  \| X_1\|^2 - 1 }^6  \IND_\cE = O (1). 
$$
To this end, for integer $\ell \geq 1$, we write
\begin{eqnarray*}
\dE\ABS{  \| X_1\|^2 - 1 }^{\ell} \IND_\cE & =  & \dE\ABS{  \| X_1\| - 1 }^{\ell} \ABS{  \| X_1\| + 1 }^{\ell} \IND_\cE  \leq  3^\ell  \dE\ABS{  \| X_1\| - 1 }^{\ell}.
\end{eqnarray*}
Then, Theorem \ref{th:GM} implies that there exists $c_\ell$ such that
$$
\dE\ABS{  \| X_1\| - 1 }^{\ell} \leq c_\ell \, p^{ - \ell/ 6}. 
$$
It follows that 
\begin{eqnarray}\label{eq:var6}
 \dE\ABS{  \| X_1\|^2 - 1 }^{\ell} \IND_\cE = O \PAR{ p^{ - \ell / 6} }. 
\end{eqnarray}
This proves \eqref{eq:T3}. It finally remains to compare $\dE \mu_D$ and $\dE \mu_{M}$. 

\subsubsection*{Step 5 : End of proof}

We set $$z_i = (  \| X_i \|^ 2  -1 ).$$ We note that for $i \ne j$, 
$$
D_{ij} = M_{ij} + \sum_{1 \leq k + \ell \leq 3}  c_{k\ell} z_i ^{k} z_j ^{\ell},
$$
for some coefficients $c_{k \ell}$ depending on $f'(2),f''(2),f'''(2)$. Note that $c_{10} = c_{01} =  f'(2)$. Similarly, 
$$
D_{ii} = M_{ii} + 2 f'(2) z_i   = M_{ii} + c_{10}  z_i + c_{01}  z_i . 
$$
 Define the matrix $E$, for all $ 1 \leq i , j \leq n$, 
$$
E_{ij} = M_{ij} +  \sum_{1 \leq k + \ell \leq 3}  c_{k\ell} z_i ^{k} z_j ^{\ell}.
$$ 
If $\cE$ holds, then $\max_i |z_i |\leq \veps(n)$ and we find
$$
\ABS{ E_{ij} - D_{ij} } = \IND ( i = j )\ABS{  \sum_{2 \leq k  + \ell \leq 3}  c_{k\ell} z_i ^{k} z_i ^{\ell} }  \leq c \IND ( i = j ) \veps(n)^2.
$$
It follows from \eqref{eq:hehe} that 
\begin{eqnarray*}
d ( \dE \mu_D , \dE \mu_E ) \leq \dE  d (  \mu_D ,  \mu_E ) & \leq & \dP ( \cE^c) + \PAR{ \frac 1 n  \sum_{i,j}  \dE  \ABS{ E_{ij} - D_{ij} }^2 \IND_\cE}^{1/2} \\
& \leq &   \dP ( \cE^c) +  c \veps(n)^2. 
\end{eqnarray*}
We deduce that
$$
\lim_{n \to \infty}  d ( \dE \mu_D , \dE \mu_E )  = 0. 
$$
We notice finally that the matrix $E - M$ is equal to 
$$
 \sum_{1 \leq k + \ell \leq 3}  c_{k\ell} {Z_k} {Z_\ell}^T, 
$$
where $Z_{k}$ is the vector with coordinates $(z_i ^{k})_{1 \leq i \leq n}$.  It implies in particular that $\rank ( E-M ) \leq 9$, indeed the $\rank$ is subadditive and  $\rank ( {Z_{k}} {Z_{\ell}}^T) \leq 1$. In particular, it follows from \eqref{eq:hehe} that 
$$
d ( \dE \mu_E , \dE \mu_M ) \leq \dE  d (  \mu_E ,  \mu_M )  \leq \frac 9 n .
$$
This concludes the proof of Proposition \ref{prop:main} and of Theorem \ref{th:main}.

\subsection{Proof of Theorem \ref{th:main2}}

The concentration inequality \eqref{eq:concineq} holds. It is thus sufficient to prove the analog of Proposition \ref{prop:main}. If $\ell \geq 2$, the proof is essentially unchanged. In step $1$, the assumption \eqref{eq:boundEE} implies the existence of a sequence $\veps = \veps(n)$ going to $0$ such that $\dP ( \cE ) \to 1$.  Then, in step $4$, it suffices to extend the Taylor expansion up to $\ell$. 

For the case $\ell =1$  : in step $2$, we perform directly the Taylor expansion around $2$, for $i \ne j$ we write $f(\| X_i - X_j \|^2) = f (2) - 2 f'(2) X_i ^T X_j ( 1 + o(1))$.  We then move directly to step $5$. (As already pointed, this case is treated in \cite{DoVu}).

\bibliographystyle{abbrv}
\bibliography{mat}

\end{document}